\documentclass[11pt]{article}

\usepackage[T1]{fontenc}
\usepackage[latin1]{inputenc}
\usepackage{geometry}
\usepackage{amsmath}
\usepackage{amssymb}
\usepackage{indentfirst}
\usepackage{pifont}
\usepackage{graphicx}
\usepackage{stmaryrd}
\usepackage[all]{xy}
\usepackage{amsthm}
\usepackage{setspace}
\usepackage{lmodern}
\usepackage{hyperref}
\newcommand{\arXiv}[1]{\href{http://arxiv.org/abs/#1}{\texttt{arXiv:#1}}}
\newcommand{\zambon}{B\cap J_1 B\cap J_2 B\cap J_1 J_2 B}

\theoremstyle{plain}
\newtheorem{thm}{Theorem}[section]
\newtheorem{lem}[thm]{Lemma}
\newtheorem{prop}[thm]{Proposition}
\newtheorem{cor}[thm]{Corollary}
\theoremstyle{definition}
\newtheorem{defn}[thm]{Definition}
\newtheorem{example}[thm]{Example}
\newtheorem{examples}[thm]{Examples}
\newtheorem{rmk}[thm]{Remark}

\DeclareMathOperator{\Ad}{Ad}
\DeclareMathOperator{\idn}{idn}

\DeclareMathOperator{\Dir}{Dir}
\DeclareMathOperator{\local}{loc}
\DeclareMathOperator{\Diffeo}{Diffeo}
\DeclareMathOperator{\REAL}{Re}

\newcommand{\rond}{\mbox{\tiny{$\circ $}}}

\newcommand{\mfgs}{\mathfrak{g}^*}
\newcommand{\mfg}{\mathfrak{g}}
\newcommand{\ip}[2]{\left< #1 , #2 \right>}
\newcommand{\Courant}[2]{\left[ #1 , #2 \right]}
\newcommand{\Lie}[2]{[ #1 , #2 ]}
\newcommand{\ld}[1]{\mathcal{L}_{#1}}

\newcommand{\cc}[1]{\overline{#1}}

\newcommand{\gendex}[2]{\left\{ #1 \right\}_{#2}}
\newcommand{\TC}{T_{\mathbb{C}}}
\newcommand{\TCS}{T^*_{\mathbb{C}}}
\newcommand{\ssi}{\Leftrightarrow}
\newcommand{\pXG}{\big( \TC M_0 \oplus (\TC \orb)^0 \big)}
\newcommand{\XG}{\TC M_0 \oplus (\TC \orb)^0}

\newcommand{\GX}{\TC \orb \oplus (\TC M_0)^0}
\newcommand{\pXGr}{\big( T M_0 \oplus (T \orb)^0 \big)}
\newcommand{\XGr}{T M_0 \oplus (T \orb)^0}
\newcommand{\pGXr}{\big( T \orb \oplus (T M_0)^0 \big)}
\newcommand{\GXr}{T \orb \oplus (T M_0)^0}
\newcommand{\backl}{\mathbin{\vrule width1.5ex height.4pt\vrule height1.5ex}}
\newcommand{\ii}{\backl}
\newcommand{\BB}{\XGr}
\newcommand{\Bp}{\GXr}
\newcommand{\bemol}{\flat}
\newcommand{\orb}{\mathcal{F}}
\newcommand{\sig}{\sigma}
\newcommand{\sigt}{\tilde{\sigma}}
\newcommand{\taut}{\tilde{\tau}}
\newcommand{\sigh}{\hat{\sigma}}
\newcommand{\secloc}[1]{\Gamma_{\local}(#1)}
\newcommand{\seclocinv}[1]{\Gamma_{\local}(#1)^G}
\newcommand{\sect}[1]{\Gamma(#1)}
\newcommand{\secinv}[1]{\Gamma(#1)^G}

\newcommand{\backward}{\mathcal{B}}
\newcommand{\forward}{\mathcal{F}}

\newcommand{\ctXG}{TM_G\oplus T^*M_G}
\newcommand{\cctXG}{\TC M_G\oplus \TC^*M_G}
\newcommand{\imaginary}{\mathbf{i}}
\newcommand{\BC}{B_{\mathbb{C}}}

\begin{document}

\title{Reduction of Generalized Complex Structures}
%\author{Mathieu Sti\'enon and Ping Xu}
\author{Mathieu Sti\'enon\thanks{Francqui fellow of the Belgian American Educational Foundation} and Ping Xu\thanks{Research partially supported by NSF grant DMS03-06665 and NSA grant 03G-142.} \bigskip \\
Department of Mathematics \\ Pennsylvania State University \\ University Park, PA 16802 \bigskip \\ \texttt{stienon@math.psu.edu} \\ \texttt{ping@math.psu.edu}}
%\\[5pt] {\small\it  Department of Mathematics, Pennsylvania State University}\\[5pt]
%{\small\it e-mail: stienon@math.psu.edu, ping@math.psu.edu}}

%\date{\texttt{file: \jobname.tex}}
\date{}
\maketitle

\begin{abstract}
We study reduction of generalized complex structures. More precisely, we investigate the following question. Let $J$ be a generalized complex structure on a manifold $M$, which admits an action of a Lie group $G$ preserving $J$. Assume that $M_0$ is a $G$-invariant smooth submanifold and the $G$-action on $M_0$ is proper and free so that $M_G:=M_0/G$ is a smooth manifold. Under what condition does $J$ descend to a generalized complex structure on $M_G$? We describe a  sufficient
 condition for the reduction to hold, which includes the Marsden-Weinstein reduction of symplectic manifolds 
and the   reduction of the complex structures in K\"ahler manifolds
 as special cases. As an application, we study reduction of
generalized K\"ahler manifolds.
\end{abstract}

\tableofcontents

\section{Introduction}

Generalized complex structures \cite{Hitchin, Gualtieri}
  have been extensively studied recently
 due to their close connection with mirror symmetry.
They include both symplectic and complex structures as extreme cases.

As it is well-known in symplectic geometry, a useful way 
of producing new symplectic manifolds is via the so called Marsden-Weinstein reduction \cite{MW}. This is a procedure for constructing symplectic manifolds from a Hamiltonian system with symmetry admitting a momentum map. Let us recall this construction briefly below. Suppose we are given a symplectic manifold $(M,\omega)$, an action of a Lie group $G$ on $M$ preserving the symplectic form, and an equivariant momentum map $\mu:M\to \mfgs$, i.e. $\mu$ satisfies the relations: 
$$\mu(g\cdot x)=\Ad^*_{g^{-1}} \mu(x), \ \forall g\in G , x\in M$$ and $$\hat{A} \ii \omega = d \ip{\mu}{A}, \quad \forall A \in \mfg ,$$ where $\hat{A}$ denotes the fundamental vector field generated by $A\in \mfg$. Assume that $0$ is a regular value of the momentum map $\mu$ so that the preimage $M_0=\mu^{-1}(0)$ is a  $G$-invariant smooth submanifold.
Moreover, if we assume 
that the $G$-action on $M_0$ is free and proper so that $M_G=M_0/G$ is a smooth manifold,
then it inherits a natural symplectic structure \cite{MW} (see \cite{SL}
for the symplectic reduction in singular case).  In the context
of Poisson manifolds, a  reduction procedure was carried out by
Marsden-Ratiu \cite{MR} and  Ortega-Ratiu in the singular 
case \cite{OR0}.  See \cite{OR} for a  beautiful comprehensive study
 on  Hamiltonian reductions.

On the other hand, since there is no such notion of momentum maps for complex manifolds, there does not exist a general scheme of reduction of $G$-invariant
complex structures in the literature as far as we know. However, for a $G$-invariant K\"ahler manifold which admits an equivariant momentum map for the symplectic structure, one can prove that the complex structure can also descend to the symplectic reduced space $M_G=M_0/G$ so that $M_G$ becomes a K\"ahler manifold. This is the so called K\"ahler reduction.

There are several versions of K\"ahler reduction, which were due to Guillemin-Sternberg \cite{GS}, Kirwan \cite{K}, and Hitchin et. al. \cite{Hitchinet.al.} respectively. A careful examination of the argument used by Guillemin-Sternberg in \cite{GS} shows that the following identity:
\begin{equation} 
T_x M=T_x M_0\oplus j(T_x (G\cdot x)) \label{3} ,\ \ \ \forall x\in M_0,
\end{equation} 
plays an essential role in carrying out the reduction for the complex structure, where $j: TM\to TM$ is the complex structure of the K\"ahler manifold. Equation \eqref{3} holds automatically when $M_0=\mu^{-1}(0)$ is the level set of the momentum map $\mu: M  \to \mfgs$. Indeed this is exactly why the symplectic reduced space  $M_G=M_0/G$ inherits a complex structure.

A natural question arises as to whether there is a reduction procedure for a $G$-invariant generalized complex structure which combines the above two special cases. More precisely, the question can be formulated as follows. 

\begin{quote}
Let $J$ be  a  generalized complex structure on a manifold $M$, which admits an action of a Lie group $G$ preserving $J$. Assume that $M_0$ is a $G$-invariant smooth submanifold and the $G$-action on $M_0$ is free and proper so that $M_G:=M_0/G$ is a smooth manifold. Under what condition does $J$ descend to a generalized complex structure on $M_G$?
\end{quote}

 In this paper, we describe a sufficient condition for such a reduction.
Our condition comprises many well-known examples as special cases including the above two important cases. In particular, even when $J$ is a honest complex structure, we derive some interesting condition for the reduction of the complex structure to hold, which seems to be new.

Below let us describe the main idea of our approach briefly. There are several equivalent definitions of generalized complex structures. A useful one for us is  that a generalized complex structure on a manifold $M$ is a pair of transversal (complex) Dirac structures $(E_+,E_-)$ on $M$ which are complex conjugate to each other. In other words $(E_+, E_-)$ constitutes a (complex) Lie bialgebroid in the sense of \cite{MX} (see \cite{LWX}). 
Roughly, our approach is as follows.
Using the inclusion map $M_0\to M$, one may pull back $E_{\pm}$ to $M_0$ to obtain a pair of complex conjugate Dirac structures on $M_0$. Since $J$ is $G$-invariant, using the push forward map $M_0\to M_0/G=M_G$, one obtains a pair of complex conjugate Dirac structures $E'_{\pm}$ on $M_G$. To ensure this gives rise to a generalized complex structure on $M_G$, a necessary  condition is
 that $E'_+\cap E'_-=0$. However, there is  a subtlety.
It is due to the difficulty that pull back and push forward
of Dirac structures may not be smooth bundles.
This forces us to impose some extra smoothness assumptions.
Indeed, in this paper we combine the pull back and push forward steps together
to derive a sufficient condition for the reduction. Note that
Dirac reduction in the real context  has been  first studied
 by Blankenstein (see \cite{Blankenstein, BS})
and by Blankenstein-Ratiu in the singular case \cite{BR}
in  connection with the study of  Hamiltonian mechanics.

 One may easily check that if one starts with a symplectic structure, the reduced generalized complex structure is still symplectic, while reduction of a complex structure is still a complex structure. Thus the two examples above are indeed special cases of our general condition.
As an application, we study  reduction of  generalized  K\"ahler  manifolds,
which generalizes the usual  K\"ahler reduction \cite{GS, Hitchinet.al., K}.
Our paper only  deals with the  analogue of regular K\"ahler reduction.
See \cite{Huebschmann} for the singular case.

The paper is organized as follows. In Section 2, we recall some basic materials in Dirac geometry. In particular, we describe the pull back and push forward constructions. In Section 3, we recall several equivalent definitions of generalized complex structures. Section 4 is devoted to the study of reduction of generalized complex structures. As corollaries, we consider several special examples including the reduction of usual complex structures, and the $B$-transform of symplectic structures. In Section 5, we give an explicit description of the reduced generalized complex structure in terms of
the endomorphism $J_G$ of the vector bundle 
$TM_G\oplus T^*M_G$. In Section 6, as an application, we investigate 
the reduction of generalized K\"ahler structures.

As a   special case when  $G=\{*\}$, our main results
would lead to conditions which guarantee submanifolds of a generalized  
complex manifold to  inherit natural generalized
complex structures. This topic was studied in detail separately by Barton and
the first author \cite{BartonS}. Note though that the first investigation of the submanifolds of a generalized complex manifold was \cite{Oren}.

The results of this paper were announced at the conference ``Poisson geometry'' held in Trieste in July, 2005, where 
we learned that Bursztyn-Cavalcanti-Gualtieri \cite{BG} were
 working on a related subject. Subsequently, we noticed that
there have appeared several papers studying similar
topics independently, including the one by Hu \cite{H}
 the one by Lin-Tolman \cite{LT} and the one by Vaisman \cite{Vaisman}. 
It would be interesting to clarify the relations between all these approaches.

\paragraph{Notations} We denote by $V_{\mathbb{C}}$ the complexified $V\otimes\mathbb{C}$ of a vector space $V$. Likewise, $\TC M$ and $\TC^* M$ respectively denote the tangent and cotangent bundles of a manifold $M$.  
If $W$ is a vector subspace of $V$, $W^0$ denotes the annihilator subspace of $W$ inside $V^*$.
And if $V$ is endowed with an inner product, $W^{\perp}$ denotes the subspace of $V$ orthogonal to $W$. Finally, $\mathbf{i}$ denotes the imaginary number $\sqrt{-1}$.

\paragraph{Acknowledgments} 
Ping Xu would like to thank
 the \emph{Universit\'e Pierre et Marie Curie, Paris} 
 for her hospitality while work on this project was being done.
 We would like to thank Henrique Bursztyn, Camille Laurent-Gengoux,
Aissa Wade and Alan Weinstein for many useful discussions. 
Special thanks go to Tudor Ratiu, who kindly called our attention to many references on related subjects we overlooked in an earlier version, and to Marco Zambon, who spotted a (minor) error.
 Finally, we thank the organizers of the conference \emph{``Poisson geometry'' - Trieste, July 2005},
 for having given us the chance to present this work.

\section{Dirac geometry}

In this section, we shortly review some basic ingredients of Dirac geometry which will be used in this paper. For details, see \cite{BursztynRadko, Courant}.

Given a vector space $V$, we consider the direct sum $V\oplus V^*$ endowed with the \emph{inner product} 
$$\ip{X+\xi}{Y+\eta}=\tfrac{1}{2}\big(\xi(Y)+\eta(X)\big), \quad \forall X,Y\in V, \; \xi, \eta \in V^* .$$
We denote by $\rho$ and $\rho^*$ the natural projections of $V\oplus V^*$ onto $V$ and $V^*$ respectively.
A \emph{Dirac structure} on $V$ is a maximal isotropic subspace of $V\oplus V^*$ with respect to $\ip{\cdot}{\cdot}$. The set of linear Dirac structures on $V$ is denoted by \emph{$\Dir(V)$}.

Alternatively, a Dirac structure on a vector space can be described as follows.
Let $L \in \Dir(V)$.
Since $L$ is isotropic with respect to $\ip{\cdot}{\cdot}$, there exists a natural skew-symmetric bilinear form $\Lambda$ on $L$ defined by $$\Lambda(X+\xi,Y+\eta)=\xi(Y)=-\eta(X), \quad \forall X+\xi, Y+\eta\in L.$$ One easily checks that 
$$\Lambda(X+\xi_1,Y+\eta_1)=\Lambda(X+\xi_2,Y+\eta_2), \quad \forall X+\xi_i, Y+\eta_i \in L, \ i=1, 2$$ and 
$$\Lambda(X_1+\xi,Y_1+\eta)=\Lambda(X_2+\xi,Y_2+\eta), \quad \forall X_i+\xi, Y_i+\eta \in L , \ i=1, 2.$$
Hence, there exists a 2-form $\Omega$ on $\rho(L)$ and a 2-form $\pi$
 on $\rho^*(L)$ defined, respectively, by 
$$\Omega(X,Y)=\Lambda(X+\xi,Y+\eta)=-\pi(\xi,\eta), \quad \forall X+\xi, Y+\eta \in L.$$

Let $X\in \rho(L)$. Then there is a $\xi\in V^*$ such that
 $X+\xi\in L$.  And we have $\Omega(X,Y)=\xi(Y), \quad \forall Y \in \rho(L)$,
 or, equivalently, $X\ii\Omega=\xi_{|\rho(L)}$,
 where $\xi_{|\rho(L)}$ denotes the restriction of $\xi\in V^*$ to $\rho(L)$. Thus 
$$ X+\xi\in L \quad \Longleftrightarrow \quad \begin{cases} X \in \rho(L) \\ X\ii\Omega=\xi_{|\rho(L)} .\end{cases} $$
Thus knowing the Dirac structure $L$ on $V$ is exactly the same as
 knowing the subspace $\rho(L)$ of $V$ and the 2-form $\Omega$  on $\rho(L)$.
Similarly, $L$ is equivalent to the pair $(\rho^*(L), \pi)$. 

Therefore, we have the converse: any subspace $R\subset V$ endowed with a 2-form $\Omega$ defines a Dirac structure 
$$\{ X+\xi | X \in R, \xi \in V^*, X\ii\Omega=\xi_{|R} \}$$ on $V$ and any subspace $R^*\subset V^*$ endowed with a 2-form $\pi$ defines a Dirac structure
$$\{ X+\xi | \xi \in R^*, X \in V, \pi(\xi,\eta)=-\eta(X) \; , \forall \eta\in R^* \}$$ on $V$.

Dirac structures can be pulled back and pushed forward using a linear
map.  Let $V\xrightarrow{\phi}W$ be a linear map.
Consider a Dirac structure $L_W$ on $W$ and its associated pair $\big(\rho(L_W),\Omega_{L_W}\big)$. The 2-form $\Omega_{L_W}$ pulls back to the preimage $\phi^{-1}(\rho(L_W))$, yielding the pair $\big(\phi^{-1}(\rho(L_W)),\phi^*\Omega_{L_W}\big)$. The Dirac structure on $V$ associated to this new pair is called \emph{backward} image of $L_W$ through $\phi$ and is denoted by $\mathcal{B}\phi(L_W)$.

Let $W^*\xrightarrow{\phi^*}V^*$ be the dual  of $\phi$. 
Consider a Dirac structure $L_V$ on $V$ and its associated pair $\big(\rho^*(L_V),\pi_{L_V}\big)$. The pullback of the 2-form $\pi_{L_V}$ to the preimage ${(\phi^*)}^{-1}\big(\rho^*(L_V)\big)$ is equal to $\phi\pi_{L_V}$. 
The Dirac structure on $W$ associated to the pair $\big({(\phi^*)}^{-1}\big(\rho^*(L_V)\big),\phi\pi_{L_V}\big)$ is called the \emph{forward} image of $L_V$ through $\phi$ and is denoted by $\mathcal{F}\phi(L_V)$.

One easily checks that 
$$\mathcal{B}\phi(L_W)=\{X+\phi^*\xi|X\in V,\xi\in W^*, \phi X+\xi\in L_W\}$$ and 
$$\mathcal{F}\phi(L_V)=\{\phi(X)+\xi|X\in V,\xi\in W^*, X+\phi^*\xi\in L_V\} .$$

Let $M$ be a smooth manifold. A \emph{Dirac structure} on $M$ is a smooth subbundle $L\subset TM\oplus T^*M$ which determines a Dirac structure in $T_m M\oplus T^*_m M$ for each $m\in M$ and whose space of sections is closed under the (skew-symmetric) \emph{Courant bracket} 
 on $\Gamma(TM\oplus T^* M)$ defined  by 
\begin{multline*} \Courant{X+\xi}{Y+\eta}=[X, Y]+\ld{X}\eta-\ld{Y}\xi+
\tfrac{1}{2}d\big(\xi(Y)-\eta(X)\big), \\ \forall X, Y
\in \mathfrak{X} (M), \ \xi, \eta \in \Omega^1 (M) .\end{multline*}
In general, the Courant bracket does not satisfy the Jacobi identity. However, it does when being
 restricted to the space $\Gamma(L)$ of sections of a Dirac subbundle $L\subset TM\oplus T^*M$.
Then $L\to M$ inherits a Lie algebroid structure with the projection $L\xrightarrow{\rho_{|L}}TM$ as anchor map and the restriction of the Courant bracket as Lie algebroid bracket 
\cite{Courant}.

Although the backward and forward images of a Dirac structure through the differential of a smooth map are always pointwise Dirac structures, they are generally not constant rank smooth vector bundles. 
The remainder of this section is devoted to the description of situations in which backward and forward images of Dirac structures are again Dirac structures.

\begin{defn}
Let $M\xrightarrow{\varphi}N$ be a smooth map.
Two generalized tangent vectors $X+\xi\in TM\oplus T^*M$ and $Y+\eta\in TN\oplus T^*N$ are said to be \emph{$\varphi$-related}, denoted by $X+\xi \overset{\varphi}{\rightsquigarrow}Y+\eta$, if $\varphi_*X=Y$ and $\varphi^*\eta=\xi$.
Clearly, this notion extends to sections.
%\item Two sections $\sig_M=X+\xi\in\sect{TM\oplus T^*M}$ and $\sig_N=Y+\eta\in\sect{TN\oplus T^*N}$ are said to be \emph{$\varphi$-related}, denoted by $\sig_M\overset{\varphi}{\rightsquigarrow}\sig_N$, if $\varphi_*X=Y$ and $\varphi^*\eta=\xi$.
%\end{enumerate}
\end{defn}

\begin{lem} \label{lem:starstar}
Assume that $\sig_M^i\in\sect{TM\oplus T^*M}$ and $\sig_N^i\in\sect{TN\oplus T^*N}$ satisfy $\sig_M^i\overset{\varphi}{\rightsquigarrow}\sig_N^i$, where $i=1,2$. Then $$\Lie{\sig_M^1}{\sig_M^2} \overset{\varphi}{\rightsquigarrow} \Lie{\sig_N^1}{\sig_N^2} ,$$ where the bracket on the LHS refers to the Courant bracket on $\sect{TM\oplus T^*M}$ while the bracket on the RHS refers to the Courant bracket on $\sect{TN\oplus T^*N}$.
\end{lem}

\begin{proof}
Write $\sig_M^i=X^i+\xi^i$ and $\sig_N^i=Y^i+\eta^i$, where $X^i\in\mathfrak{X}(M)$, $\xi^i\in\Omega^1(M)$, $Y^i\in\mathfrak{X}(N)$ and $\eta^i\in\Omega^1(N)$ for $i=1,2$. Since
$\sig_M^i$ and $\sig_N^i$
  are $\varphi$-related, $\varphi_*X^i=Y^i$ and $\varphi^*\eta^i=\xi^i$ for $i=1,2$.
Then $$\Lie{\sig_M^1}{\sig_M^2}=\Lie{X^1}{X^2}+\ld{X^1}\xi^2-\ld{X^2}\xi^1+\tfrac{1}{2}d\big(\xi^1(X^2)-\xi^2(X^1)\big)$$ and $$\Lie{\sig_N^1}{\sig_N^2}=\Lie{Y^1}{Y^2}+\ld{Y^1} \eta^2-\ld{Y^2}\eta^1+\tfrac{1}{2}d\big(\eta^1(Y^2)-\eta^2(Y^1)\big) .$$
Now  $$\varphi_*\Lie{X^1}{X^2}=\Lie{\varphi_*X^1}{\varphi_*X^2}=\Lie{Y^1}{Y^2} ,$$ 
and
\begin{align*} \varphi^*(\ld{Y^1}{\eta^2}) = & \varphi^*(i_{Y^1}d\eta^2+di_{Y^1}\eta^2) \\ 
=& \varphi^* i_{\varphi_* X^1}d\eta^2+\varphi^* d\big(\eta^2(Y^1)\big) \\ 
=& i_{X^1}\varphi^* d\eta^2+d\big(\xi^2(X^1)\big) \\ 
=& i_{X^1}d\xi^2+d i_{X^1}\xi^2 \\ 
=& \ld{X^1}\xi^2 .\end{align*}
Similarly, we have 
$$\varphi^*(\ld{Y^2}{\eta^1}) = \ld{X^2}\xi^1 .$$
Finally 
$$\varphi^*d\big(\eta^1(Y^2)-\eta^2(Y^1)\big)=d\big(\xi^1(X^2)-\xi^2(X^1)\big)$$ since $$\varphi^*\big(\eta^1(Y^2)\big)=\varphi^*\big(\eta^1(\varphi_*X^2)\big)=(\varphi^*\eta^1)(X^2)=\xi^1(X^2) .$$
This concludes the proof.
\end{proof}

Recall that a \emph{generalized smooth subbundle} of a vector bundle
$E\to M$ is a subbundle $V\subset E$ (whose fibers may not be of constant ranks)
such that  for any vector $v\in V_m$ there is a   a local smooth section
$s$ of $V$ extending $v$, i.e. $s|_m=v$.

\begin{prop}
\label{pro:2.3}
\begin{enumerate}
\item Let $M\xrightarrow{i}N$ be an embedding and $L_N$ a Dirac structure on $N$. If, for any vector $v_m\in\backward i(L_N)$, there exists a local smooth section $s$ of $L_N$ defined in a neighborhood of $i(m)\in N$, such that $s_{|M}\in \secloc{i_*TM\oplus T^*N|_M}$ and $v_m\overset{i}{\rightsquigarrow}s|_{i(m)}$, then $\backward\phi(L_N)$ is a Dirac structure on $M$.
In particular, the above condition holds if 
$L_N|_M\cap\big(i_*(TM)\oplus T^*N|_M\big)$ is a
generalized smooth subbundle.
\item Let $M$ be a smooth manifold endowed with a proper and free action of a Lie group $G$. Let $\pi$ denote the canonical projection of $M$ onto $N:=M/G$. If $L_M$ is a $G$-invariant Dirac structure on $M$, and for any vector $v_n\in\forward\pi(L_M)$, there exists a local smooth section $s$ of $L_M$ defined in a neighborhood of a point $m\in M$ in the fiber of $\pi$ over $n$ such that $s\in\secloc{TM\oplus\pi^*(T^*N)}$ and $s|_m\overset{\pi}{\rightsquigarrow}v_n$, for instance when
$L_M\cap\big(TM\oplus\pi^*(T^*N)\big)$ is a generalized smooth subbundle, then $\mathcal{F}\pi(L_M)$ is a Dirac structure on $N$.
\end{enumerate}
\end{prop}
\begin{proof} 
(1). First, we observe that $\mathcal{B}i (L_N)$ is well-defined
 since $i$ is injective. The assumption implies that
any vector $v_m \in \backward i (L_N)$ admits a local
section $\tilde{s}$  of $\backward i (L_N)$
such that $\tilde{s}\overset{i}{\rightsquigarrow} s$, where
$s$ is  a local smooth section of $L_N$ defined in
a neighborhood of $i (m)\in N$ such that $s_{|M}\in \secloc{i_*TM\oplus T^* N|_M}$. 
Then $\tilde{s}$ must be a smooth
local section. Therefore  $\backward i (L_N)$ is a smooth
subbundle. 

Now for any point $m\in M$, take a basis $\gendex{v^k\in \backward i (L_N)|_m }{k\in K}$.
Let  $\tilde{s}^k$ be a local section  of $\backward i (L_N)$
through $v^k$ as above. Then $\tilde{s}^k\overset{i}{\rightsquigarrow} s^k$,
where $s^k$ is a local smooth section of $L_N$. Thus $[\tilde{s}^k, \tilde{s}^l]$, $\forall k, l$, is still a local smooth section  of $\backward i (L_N)$
according to Lemma  \ref{lem:starstar}.
Hence it follows that  $\backward i (L_N)$ is Courant involutive.

For the second part, note that any element $v_m\in \backward i (L_N)$ naturally
induces an element $v_{i(m)}\in  L_N|_M\cap\big(i_*(TM)\oplus T^*N|_M\big)$
such that $v_m \overset{i}{\rightsquigarrow} v_{i(m)}$.
Thus the condition is satisfied since $L_N|_M\cap\big(i_*(TM)\oplus T^*N|_M\big)$ is a generalized smooth subbundle by assumption.

% Fix $m\in M$. Since $L_N\cap\big(i_*(TM)\oplus T^*N\big)$ is weakly smooth
%  and $\mathcal{B}i(L_N)$ has constant rank, there exists a finite set 
%  $\gendex{i_*X_k+\xi_k}{k\in K}$ of local sections of $L_N\cap\big(i_*(TM)\oplus T^*N\big)$ such that the linear combinations of elements in $\gendex{X_k+i^*\xi_k}{k\in K}$ with coefficients in $C^{\infty}(M)$ generate $\sect{\mathcal{B}i(L_N)}$ locally around $m$. Hence $\mathcal{B}i(L_N)$ is a smooth vector bundle.

(2). First, we observe that $\mathcal{F}\pi(L_M)$ is well-defined
 since $L_M$ is $G$-invariant. Now
 for any vector $v_n\in\forward\pi(L_M)$,
let  $s$ be  a local section of $L_M$ defined in a neighborhood
 of a point $m\in \pi^{-1}(n)$ such that 
$s\in \secloc{TM\oplus\pi^*(T^*N)}$ and $s|_m 
\overset{\pi}{\rightsquigarrow}v_n$. 
Let $L\subset M$ be
a slice through $m$ transversal to the $G$-orbits. 
Since $L_M$ is  $G$-invariant, there is a local
smooth section $\hat{s}$ of $L_M$ such that
$\hat{s}|_L=s|_L$. It is clear that $
\hat{s}\in \secloc{TM\oplus\pi^*(T^*N)}$. Hence  it induces
a local section $\tilde{s}$ of $\forward\pi(L_M)$ through
$v_n$. Therefore $\forward\pi(L_M)$ is indeed
a smooth subbundle.  The rest can be proved similarly as in (1).
\end{proof}

\section{Generalized complex structures}
\label{sec:gcs}

Let $V$ be a vector space. Consider a linear endomorphism
 $J$ of $V\oplus V^*$ such that $J^2=-I$ and $J$ is orthogonal with respect to the inner product $$\ip{X+\xi}{Y+\eta}=\tfrac{1}{2}\big(\xi(Y)+\eta(X)\big), \quad \forall X,Y\in V, \; \xi, \eta \in V^* .$$ Such a linear map is called a \emph{linear generalized complex structure} by Hitchin and Gualtieri \cite{Gualtieri, Hitchin}. The complexified vector space $(V\oplus V^*)\otimes \mathbb{C}$ decomposes 
as the direct sum $$(V\oplus V^* )\otimes \mathbb{C}=E_+\oplus E_-$$ 
of the eigenbundles of $J$ corresponding to the eigenvalues $\pm \imaginary$
 respectively, i.e.,
$$E_{\pm}=\{(X+\xi)\mp \imaginary J(X+\xi)|X+\xi\in V\oplus V^* \} .$$ Both eigenspaces are 
maximal isotropic with respect to
 $\ip{\cdot}{\cdot}$ and they are complex conjugate to each other.

The following lemma is obvious.
\begin{lem}
The linear generalized complex structures are in 1-1 correspondence with the splittings 
$(V\oplus V^*)\otimes \mathbb{C}=E_+\oplus E_-$ with $E_{\pm}$ maximal isotropic and $E_-=\cc{E_+}$.
\end{lem}

Now, let $M$ be a manifold and $J$ a bundle endomorphism
 of $TM\oplus T^*M$ such that $J^2=-I$,  and $J$ is 
orthogonal with respect to  $\ip{\cdot}{\cdot}$. 
In the associated eigenbundle decomposition $$(T_\mathbb{C} M\oplus T_\mathbb{C}^*M)
=E_+\oplus E_- ,$$ if $\Gamma(E_+)$ is closed under the (complexified) Courant bracket, then $E_+$ is a (complex) Dirac structure on $M$ and one says that $J$ is a 
\emph{generalized complex structure}\cite{Gualtieri, Hitchin}.
In this case, $E_-$ must also be  a Dirac structure
since $E_-=\cc{E_+}$. Indeed $(E_+, E_-)$ is a complex Lie bialgebroid in the sense  of
\cite{MX}, in which $E_+$ and $E_-$ are complex conjugate to each other. 

The following proposition gives two equivalent definitions of
a generalized complex structure.

\begin{prop}
A generalized complex structure is equivalent to any of the following:
\begin{enumerate}
\item A bundle endomorphism  $J$ of  $TM\oplus T^*M$ such that $J$ is
orthogonal with respect to $\ip{\cdot}{\cdot}$ and  $J^2=-I$
satisfying
\begin{equation} \Lie{Je_1}{Je_2}-\Lie{e_1}{e_2}-J(\Lie{Je_1}{e_2}+\Lie{e_1}{Je_2})=0, \quad \forall e_1, e_2\in
\Gamma(TM\oplus T^*M) . \label{nijenhuis} \end{equation}
\item A complex Lie bialgebroid $(E_+, E_-)$ such that its
double is the standard Courant algebroid  $\TC M\oplus \TC^*M$, and $E_+$ and $E_-$ are complex conjugate to each other.
\end{enumerate}
\end{prop}
\begin{proof}
\begin{description}
\item[1$\implies$2] As above, $E_{\pm}$ are the $(\pm \imaginary)$-eigenbundles of $J$. 
Since Eq. \eqref{nijenhuis} exactly means that 
$$\Lie{e_1\mp \imaginary Je_1}{e_2\mp \imaginary Je_2}\in \Gamma(E_\pm),\quad \forall e_1, e_2 \in \Gamma(TM\oplus T^*M) ,$$
it is equivalent to the involutivity of $\Gamma(E_\pm)$.
\item[2$\implies$1] Define $J:T_\mathbb{C} M\oplus T_\mathbb{C}^*M
\to T_\mathbb{C} M\oplus T_\mathbb{C}^*M$ by setting
$J|_{E_+}=\imaginary I$ and $J|_{E_-}=-\imaginary I$. Since $E_+$ and $E_-$ are
complex conjugate to each other, we have $J\bar{e}=\overline{Je}, \ \forall e\in
T_\mathbb{C} M\oplus T_\mathbb{C}^*M$. Therefore $J$ must be
the $\mathbb{C}$-linear extension of an endomorphism
of $TM\oplus T^*M$. It is easy to check that $J$ indeed
satisfies all the axioms as in (1).
\end{description}
\end{proof}

The following are two standard examples \cite{Gualtieri, Hitchin}.
\begin{examples}
\begin{enumerate}
\item Let $j$ be an almost complex structure on $M$.
 Then $$J=\begin{pmatrix} j & 0 \\ 0 & -j^* \end{pmatrix}$$ is 
$\ip{\cdot}{\cdot}$-orthogonal and satisfies $J^2=-I$.
$J$ is a generalized complex structure if and only if
 $j$ is integrable.
\item Let $\omega$ be a nondegenerate  2-form on $M$.
 Then $$J=\begin{pmatrix} 0 & -\omega_{\bemol}^{-1} \\ \omega_{\bemol}
 & 0 \end{pmatrix}, $$ 
where $\omega_{\bemol} : TM\to T^*M$ is the bundle map $X\mapsto X\ii \omega$, is a generalized complex structure if and only if $d\omega=0$, i.e., $\omega$ is a symplectic $2$-form.
\end{enumerate}
\end{examples}

\section{Reduction using Dirac structures}
\label{sec:rds}

Let $M$ be a manifold, $J$ a generalized complex structure on $M$
 and $\varphi$ an action of a Lie group $G$ on  $
M$ preserving the generalized complex structure $J$. 
In other words, we have a group homomorphism $\varphi:G\to\Diffeo(M):g\mapsto\varphi_g$ such that, $\forall g\in G$, the bundle endomorphism
%$(\phi_g=\varphi_{g*}\oplus\varphi^*_{g^{-1}}, \varphi_g)$ of $TM\oplus T^* M\to M$ 
$$\xymatrix{
TM\oplus T^*M \ar[r]^{\phi_g} \ar[d] & TM\oplus T^*M \ar[d] \\ 
X \ar[r]^{\varphi_g}  & X }$$ (defining the induced $G$-action $\phi$ on $TM\oplus T^*M$) 
commutes with $J$, where $\phi_g=\varphi_{g*}\oplus\varphi^*_{g^{-1}}$.
 As in Section \ref{sec:gcs}, by $E_{\pm}$ we  denote  the $(\pm \imaginary)$-eigenbundles of $J$. 

Let $M_0$ be a $G$-invariant submanifold of $M$ with a free and proper $G$-action 
so that the quotient space $M_G:= M_0/G$ is a smooth manifold. Thus one
has the following maps:
$$M\overset{i}{\hookleftarrow} M_0 \overset{\pi}{\twoheadrightarrow} M_G$$

The main question that we investigate in this section is: under
what condition does $M_G$ inherit a generalized complex structure ? 

Since  $\pi$ is  surjective and $J$ is $G$-invariant,
  $E'_{\pm}:=\mathcal{F}\pi\big(\mathcal{B}i(E_{\pm})\big)$ are {well-defined, complex conjugate pointwise Dirac structures} on $M_G$. It is simple to check
that
$$E'_{\pm} =\{\pi_*X+\xi'|X\in T_{\mathbb{C}}M_0, \xi'\in T^*_{\mathbb{C}} M_G \text{ and } \exists \xi\in T^*_{\mathbb{C}}M \text{ st } i_*X+\xi\in E_{\pm} \text{ and } \pi^*\xi'=i^*\xi \} .$$ 

One possible way to ensure that $E'_{\pm}$ are Dirac structures is to use Proposition \ref{pro:2.3} by considering the pull back Dirac structure $\mathcal{B}i(E_{\pm})$ and then the push forward structure $\mathcal{F}\pi\big(\mathcal{B}i(E_{\pm})\big)$.
Below, however, we will combine these two steps together,
which enables us to obtain a stronger result.

%Assume momentaneously that $E'_+$ and $E'_-$ are Courant involutive smooth vector bundles. If $T_{\mathbb{C}}M_G\oplus T^*_{\mathbb{C}}M_G=E'_+\oplus E'_-$, one can think of them as the $(+i)$- and $(-i)$-eigenbundles of a new "reduced" generalized complex structure $J_G$ on $M_G$. 

In what follows, the foliation of $M_0$ defined by the $G$-orbits will be denoted by $\orb$, and $(T\orb)^0$ and $(TM_0)^0$ will denote the annihilators of $T\orb$ and $TM_0$ in $TM$ along the submanifold $M_0$.

It is clear that any element $i_*X_m+\xi_{i(m)}\in TM_0\oplus(T\orb)^0$ induces an element $\pi_*X_m+\xi'_{\pi(m)}\in TM_G\oplus T^*M_G$ where $\xi'_{\pi(m)}$ is defined by $i^*\xi_{\pi(m)}=\pi^*\xi'_{\pi(m)}$. Indeed, one has a bundle map 
\begin{equation} \xymatrix{
TM_0\oplus (T\orb)^0 \ar[r]^{\Pi} \ar[d] & TM_G\oplus T^*M_G \ar[d] \\
M_0 \ar[r]_{\pi}  & M_G } \label{eq:liftbundlemap} .\end{equation}
The element $i_*X_m+\xi_{i(m)}\in TM_0\oplus (T\orb)^0$ is called a \emph{lift} of $\pi_* X_m+\xi'_{\pi(m)}\in TM_G\oplus T^*M_G$. Note that, for a given element $v_n\in TM_G\oplus T^*M_G$, there exists many different lifts. Indeed, there is a choice of point $m$ in the fiber $\pi^{-1}(n)\subset M_0$ involved and $\Pi$ has $T\orb\oplus (TM_0)^0$ as non trivial kernel.

\begin{prop} \label{prop:condition}
The following assertions are equivalent.
\begin{enumerate}
\item \label{enumi} $E'_+\cap E'_-=0$
\item \label{enumii} $\big(E_++T_{\mathbb{C}}\orb\big)\cap\big(T_{\mathbb{C}}M_0\oplus (T_{\mathbb{C}}\orb)^0\big)\cap \big(E_-+(T_{\mathbb{C}}M_0)^0\big)\subset T_{\mathbb{C}}\orb\oplus (T_{\mathbb{C}}M_0)^0$
\item \label{enumiii} If $z\in T_{\mathbb{C}}\orb\oplus(T_{\mathbb{C}}M_0)^0$
 with $z=z_+ +z_-$ such that $z_{\pm}\in E_{\pm}$ and $z_{\pm}\in T_{\mathbb{C}}M_0\oplus (T_{\mathbb{C}}\orb)^0$, then $z_{\pm}\in T_{\mathbb{C}}\orb \oplus (T_{\mathbb{C}}M_0)^0$.
\item \label{enumiv} $J\big(T\orb\oplus(T M_0)^0\big)\cap 
\big(T M_0\oplus (T\orb)^0\big) \subset T\orb\oplus (T M_0)^0$
\item \label{enumv} $J\pXGr\subset\pXGr + J\pGXr$
\item \label{enumvi} $\XGr=\pXGr\cap J\pXGr + \pGXr$
\item \label{enumvii} $\TC M_G\oplus \TC^*M_G =E'_++E'_-$
\end{enumerate}
\end{prop}

\begin{proof}
\begin{description}
\item[\ref{enumi} $\ssi$ \ref{enumii}]
Note that 
\begin{multline*} E'_+\cap E'_-=\{X'+\xi'\in \TC M_G \oplus \TCS M_G | \exists X_{\pm}\in T_{\mathbb{C}} M_0 \text{ and } \xi_{\pm}\in T^*_{\mathbb{C}}M_{|M_0} \\ \text{ such that } i_* X_{\pm} + \xi_{\pm}\in E_{\pm}, \; 
\pi_* X_{\pm} =X' \text{ and } i^*\xi_{\pm} = \pi^*\xi' \} . \end{multline*}
Thus, in particular, we have 
\begin{align*} & X_{\pm}\in \TC M_0 & & \xi_{\pm}\in(\TC \orb)^0 \\ 
& X_--X_+\in \TC \orb & & \xi_--\xi_+\in (\TC M_0)^0 .\end{align*}
Assume that there exists an element $X'+\xi'\ne 0$ in $E'_+\cap E'_-$. Then there exist $i_* X_{\pm}+\xi_{\pm} \in E_{\pm}$ with $X_{\pm}\in T_{\mathbb{C}} M_0$ and $\xi_{\pm}\in T^*_{\mathbb{C}}M_{|M_0}$ as above.
% which both belong to $\TC M_0 \oplus (\TC \orb)^0$.
In particular, 
 $$i_* X_{-}+\xi_{+} \in  \pXG .$$ 
Since $\TC \orb\subset \TC M_0$, we get 
$$i_* X_- + \xi_+ = (i_* X_+ + \xi_+) + i_*(X_- -X_+) \in  E_+ + \TC \orb .$$ 
On the other hand, 
%since $(\TC M_0)^0\subset (\TC \orb)^0$, we get 
we have
$$i_* X_- + \xi_+ = (i_*X_- + \xi_-)+(\xi_+-\xi_-) \in  E_- + (\TC M_0)^0  .$$
Hence $$i_* X_- + \xi_+ \in \big( E_+ + \TC \orb\big) \cap \pXG \cap \big( E_- + (\TC M_0)^0 \big) .$$
Now $X'+\xi'\ne 0$ implies that $i_*X_-+\xi_+$ must not belong to $\GX$. 
Thus, we have proved that $E'_+\cap E'_-\ne 0$ implies 
$$\big(E_++T_{\mathbb{C}}\orb\big)\cap \pXG \cap \big(E_-+(T_{\mathbb{C}}M_0)^0\big) \nsubseteq \GX .$$ It is an easy matter to show that the converse is also true.
\item[\ref{enumii} $\Rightarrow$ \ref{enumiii}]
Take $\hat{A}+f\in \GX$ 
 and assume that  $\hat{A}+f=z_+ +z_-$ with $z_{\pm}\in E_{\pm}$
 and $z_{\pm}\in \XG$. Then $z_+-\hat{A}=-z_-+f$ is
 an element of the triple intersection of the LHS of (\ref{enumii}). Hence $z_{\pm}\in \GX$.
\item[\ref{enumiii} $\Rightarrow$ \ref{enumii}]
Any element $z$ in the triple intersection can be written both as $v+\hat{A}$ and $w+f$ for some $v\in E_+$, $\hat{A}\in \TC \orb$, $w\in E_-$ and $f\in (\TC M_0)^0$. We have $w-v=\hat{A}-f\in\GX$. The assumption (\ref{enumiii})
 implies that $v,w\in \GX$. And, finally, $z=v+\hat{A}\in\GX$.
\item[\ref{enumiii} $\ssi$ \ref{enumiv}]
Since $T_{\mathbb{C}}M=E_+\oplus E_-$, the only solution to $z=z_+ +z_-$,
where $z\in \GX$ and $z_{\pm}\in E_{\pm}$, is $z_{\pm}=\tfrac{1}{2} (z\mp \imaginary Jz)$. Since $\GX\subset\XG$, one has $z_{\pm}\in\XG$ if and only if $Jz\in \XG$. 
Thus condition (\ref{enumiii}) can be rewritten as 
$$ \left. \begin{array}{c} z\in \GX \\ Jz\in \XG \end{array} \right\} \; \Longrightarrow \; Jz \in \GX $$
which is equivalent to (\ref{enumiv}).
\item[\ref{enumiv} $\ssi$ \ref{enumv}]
Dualising one condition by taking
the annihilator with respect to $\ip{\cdot}{\cdot}$ and applying $J$ to both sides yields the other one.
\item[\ref{enumv} $\ssi$ \ref{enumvi}] Since $J^2=-I$, (\ref{enumv})
 is equivalent to $$\XGr\subset J\pXGr+\pGXr .$$ 
Remembering that $\GXr\subset\XGr$, the equivalence with (\ref{enumvi}) is obvious.
\item[\ref{enumvi} $\Rightarrow$ \ref{enumvii}] 
Any element $v\in\TC M_G\oplus \TC^*M_G$ has a lift $\tilde{v}\in\XG$.
 By (\ref{enumvi}), there
 exists $k\in\GX$ such that $\tilde{v}-k\in\pXG\cap J\pXG$.
 Then the vectors $\tfrac{1}{2}(I-\imaginary J)(\tilde{v}-k)$ and
 $\tfrac{1}{2}(I+\imaginary J)(\tilde{v}-k)$ belong to $E_+\cap\pXG$ and $E_-\cap\pXG$
 respectively and are thus the lifts of a pair of
 vectors $v_+\in E'_+$ and $v_-\in E'_-$ such that $v=v_++v_-$.
\item[\ref{enumvii} $\Rightarrow$ \ref{enumvi}]
Let $\hat{v}$ be a vector of $\XGr$. It is a lift of some 
vector $v\in TM_G\oplus T^*M_G$. Hence $v$ writes as a sum $v=v_++v_-$,
 where $v_+\in E'_+$ and $v_-\in E'_-$. Choose lifts $\tilde{v}_+\in E_+
\cap\pXG$ and $\tilde{v}_-\in E_-\cap\pXG$ of $v_+$ and $v_-$ respectively.
Then $\tilde{v}_++\tilde{v}_-$ is a lift of $v$ and belongs to
$\pXG\cap J\pXG$. And $\tilde{v}= \text{Re}(\tilde{v}_++\tilde{v}_-)$ is still
 a lift of $v$ since $v$ is a real vector. Moreover $\tilde{v}
\in \pXGr\cap J\pXGr$.
Let $k=\hat{v}-\tilde{v}$. Then it follows that
$k\in \GXr$ since both $\hat{v}$ and $\tilde{v}$ are lifts of $v$.
\end{description}
\end{proof}

As an immediate consequence, we have

\begin{cor} \label{cor:directsum}
The direct sum decomposition $\TC M_G\oplus \TC^*M_G=E'_+\oplus E'_-$ holds if and only if any element of $TM_G\oplus T^*M_G$ has a lift in $\pXGr\cap J\pXGr$. 
\end{cor}
\begin{proof}
By Proposition \ref{prop:condition}, 
the direct sum decomposition holds if and only if 
%$E'_+\cap E'_-=0$ and $\TC M_G\oplus\TC^*M_G=E'_++E'_-$.
% By proposition \ref{prop:condition}, this is equivalent to
 $$\XGr=\pXGr\cap J\pXGr + \pGXr .$$
 The result thus follows from the fact that
 $\ker\Pi=\GXr$ and the surjectivity of $\Pi$, where $\Pi$ is the map as in Eq. \eqref{eq:liftbundlemap}.
\end{proof}

It is clear that any $G$-invariant (local) section $\sigt$ of $\BB\to M_0$ induces a (local) section $\sig$ of $TM_G\oplus T^*M_G\to M_G$. Here, $G$-invariant local sections must be understood in the following obvious sense: $\sigt(m_1)=\phi_g\cdot\sigt(m_2)$, for all $m_1,m_2\in U$ such that $m_1=g\cdot m_2$ for some $g\in G$.
Hence, the bundle map \eqref{eq:liftbundlemap} induces a map $$\seclocinv{\BB}\longrightarrow\secloc{TM_G\oplus T^*M_G} .$$

\begin{defn}
%\begin{enumerate}
%\item
%A section $\sigt\in\secinv{TM\oplus T^*M}$ is said to be a 
%\emph{lift of a section} $\sig\in\sect{TM_G\oplus T^*M_G}$ 
%if $\sigt_{|M_0}\in\secinv{\BB}$ and 
%the diagram 
%\begin{equation}
%\label{eq:Pi}
%\xymatrix{
%\BB \ar[r]^{\Pi} & TM_G\oplus T^*M_G \\ 
%M_0 \ar[u]^{\sigt_{|M_0}} \ar[r]_{\pi}  & M_G \ar[u]_{\sig} %}
%\end{equation}
%commutes.
%\item 
Let $U$ be an open neighborhood in $M$. A local section $\sigt\in\secinv{(TM\oplus T^*M)_{|U}}$ is said to be a \emph{lift of a local section} $\sig\in\sect{(TM_G\oplus T^*M_G)_{|\pi(U\cap M_0)}}$ if $\sigt_{|U\cap M_0}\in\sect{(\BB)_{|U\cap M_0}}$ and the diagram 
\begin{equation}\label{eq:Pi}
\xymatrix{
\pXGr_{|U\cap M_0} \ar[r]^{\Pi} & (TM_G\oplus T^*M_G)_{|\pi(U\cap M_0)} \\ 
U\cap M_0 \ar[u]^{\sigt_{|U\cap M_0}} \ar[r]_{\pi} & \pi(U\cap M_0) \ar[u]_{\sig} }
\end{equation}
commutes.
%\end{enumerate}
\end{defn}

\begin{lem} \label{lem:liftbracket} 
Given any local sections $\sig_1, \sig_2\in\secloc{TM_G\oplus T^*M_G}$, assume that $\sigt_1, \sigt_2\in\seclocinv{TM\oplus T^*M}$ are any of their lifts. Then 
\begin{enumerate}
\item $\Lie{\sigt_1}{\sigt_2}\in\seclocinv{TM\oplus T^*M}$,
\item $\Lie{\sigt_1}{\sigt_2}_{|M_0}\in\seclocinv{\BB}$,
\item $\Pi\rond\Lie{\sigt_1}{\sigt_2}_{|M_0}=\Lie{\sigma_1}{\sigma_2}\rond\pi$, 
\end{enumerate}
where the Courant brackets are taken in $\secloc{TM\oplus T^*M}$ and $\secloc{TM_G\oplus T^*M_G}$ as is appropriate. In other words, $\Lie{\sigt_1}{\sigt_2}$ is a lift of $\Lie{\sig_1}{\sig_2}$.
\end{lem}

\begin{proof}
\begin{enumerate}
\item Assume that $\sigt_1=X+\xi$ and $\sigt_2=Y+\eta$, where $X,Y\in \mathfrak{X}_{\local}(M)^G$ and $\xi,\eta\in\Omega^1_{\local}(M)^G$. Then $$\Lie{\sigt_1}{\sigt_2}=\Lie{X}{Y}+\ld{X}\eta-\ld{Y}\xi+\tfrac{1}{2}d\big(\xi(Y)-\eta(X)\big)$$ is clearly $G$-invariant.
\item First $\Lie{X}{Y}_{|M_0}\in\mathfrak{X}(M_0)$ since $X$ and $Y$ belong to $\mathfrak{X}(M_0)$. Now, $$(\ld{X}\eta)(\hat{A})_{|M_0}=X\big(\eta(\hat{A})\big)_{|M_0}-\eta(\Lie{X}{\hat{A}})_{|M_0}, \quad \forall A\in\mfg .$$ Since $\eta(\hat{A})_{|M_0}=0$ and $X$ is tangent to $M_0$, we have $X\big(\eta(\hat{A})\big)_{|M_0}=0$. Moreover, $\Lie{X}{\hat{A}}=0$ because $X$ is $G$-invariant. Hence $(\ld{X}\eta)(\hat{A})_{|M_0}=0$. In other words, $\ld{X}\eta_{|M_0}\in (T\orb)^0$. Similarly, $\ld{Y}\xi_{|M_0}\in (T\orb)^0$. Finally, $d\big(\xi(Y)-\eta(X)\big)\in (T\orb)^0$. For $\xi(Y)-\eta(X)$ is a $G$-invariant function. Therefore we have $\Lie{\sigt_1}{\sigt_2}_{|M_0}\in \seclocinv{\BB}$.
\item By $\sigh_1, \sigh_2$, we denote the sections of $\secloc{TM_0\oplus T^*M_0}$ obtained by considering  the cotangent parts 
of $\sigt_1$ and $\sigt_2$ as elements in $T^* M_0$.
 It is clear that $\sigt_i\overset{i}{\rightsquigarrow}\sigh_i$ and $\sigh_i\overset{\pi}{\rightsquigarrow}\sig_i$. The conclusion follows from Lemma \ref{lem:starstar}.
\end{enumerate}
\end{proof}

\begin{prop} \label{prop:smoothness} 
Assume that, given any element $v_n\in TM_G\oplus T^*M_G$, there exists a
 local smooth  section $\zeta$ of $\pXGr\cap J\pXGr\to M_0$ around some
point $m \in \pi^{-1}(n)$ such that $\zeta|_m$ is a lift of $v_n$. Then the following assertions hold:
\begin{enumerate}
\item The subbundles $E'_{\pm}\subset \TC M_G\oplus\TC^*M_G$ are smooth.
\item Any smooth (global) section of $\sect{E'_{\pm}}$ locally admits a lift to a local section in $\secinv{E_{\pm}}$.
\end{enumerate}
\end{prop}

\begin{proof}
\begin{enumerate}
\item \label{enu1} Since $E'_{\pm}$ has constant rank, it suffices to prove that any element $v_n\in E'_{\pm}$ admits a smooth local section passing through it. Let $\zeta$ be a local section of $\pXGr\cap J\pXGr$ defined in an open neighborhood $V$ containing $m$ with $\pi(m)=n$ such that $\zeta|_m$ is a lift of $v_n$. Shrinking $V$ if necessary, we can extend $\zeta$ to a local section of $TM\oplus T^*M\to M$, denoted by the same symbol $\zeta$, and defined in an open neighborhood $U$ of $M_0$ such that $U\cap M_0=V$. Now take a slice $L$ in $U$ through the point $m$ and transversal to the $G$-orbits. Such a slice always exists if $U$ is taken small enough. Let $\sigt$ be the unique $G$-invariant section of $\sect{(TM\oplus T^*M)_{|U}}$ such that $\sigt_{|L}=\zeta_{|L}$. Set $\taut=\tfrac{1}{2}(I\mp \imaginary J)\rond\sigt$. 
It is clear that $\taut\in\sect{E_{\pm|U}}$ is a local $G$-invariant smooth section such that $\taut_{|U\cap M_0}\in\seclocinv{\BB}$. Therefore, it induces a smooth section $\tau\in\sect{E'_{\pm|\pi(U)}}$ 
through the element $v_n$. 
\item Let $\nu\in \sect{E'_+}$ be any smooth
section around the point $n\in M_G$. Choose a basis $\gendex{v^i_n}{i=1,\dots,k}$ of $E'_{+|n}$. Let $\tau^i$ be the local smooth section of $E'_+$ through the element $v^i_n$ constructed as in (\ref{enu1}). Then $\gendex{\tau^i|_{n'}}{i=1,\dots,k}$ is a basis of $E'_{+|n'}$ when $n'$ is in a sufficiently small neighborhood $V$ around $n$. Thus $\nu=\sum f_i \tau^i$ for some $f_i\in C^{\infty}(V)$. Let $\taut^i\in\Gamma(E_{\pm}|_U)$ be the lift of $\tau^i$ as in (\ref{enu1}) defined on an open neighborhood $U$ such that $U\cap M_0=\pi^{-1}(V)$. Choose $\tilde{f}_i\in C^{\infty}(U)^G$ such that $\tilde{f}_{i|M_0\cap U}=\pi^*f_i$. Such $\tilde{f}_i$ always exist if $U$ is chosen sufficiently small. Then $\tilde{\nu}=\sum\tilde{f}_i\taut^i$ is a lift of $\nu$ in $\secinv{E_{\pm}}$.
\end{enumerate}
\end{proof}

We are now ready to state the main theorem of this paper.

\begin{thm} \label{thm:main}
Let $M$ be a manifold endowed with an action of a Lie group $G$. Let $J$ be a $G$-invariant generalized complex structure on $M$. And let $M_0$ be a $G$-invariant submanifold of $M$ where the $G$-action is free and proper so that the quotient $M_G=M_0/G$ is a manifold. Assume that, given any element $v_n\in TM_G\oplus T^*M_G$, there exists a local smooth section $\zeta$ of $\pXGr\cap J\pXGr\to M_0$ around some point $m\in\pi^{-1}(n)$ such that $\zeta|_m$ is a lift of $v_n$. Then there is an induced generalized complex structure $J_G$ on $M_G$.
\end{thm}
\begin{proof}
This follows from Corollary~\ref{cor:directsum}, Proposition~\ref{prop:smoothness} and Lemma~\ref{lem:liftbracket}.
\end{proof}

\begin{lem} \label{prop:wsmooth}
The hypothesis of Proposition \ref{prop:smoothness} is satisfied when 
\begin{equation}\label{5bis} 
\pXGr\cap J\pXGr \text{is a generalized smooth subbundle} 
\end{equation} 
and
\begin{equation}\label{4} 
J\pGXr\cap\pXGr\subset\GXr 
.\end{equation}
\end{lem}
\begin{proof}
For any vector $v_n\in TM_G\oplus T^*M_G$, according to Proposition~\ref{prop:condition} and Corollary~\ref{cor:directsum}, Eq.~\eqref{4} implies the existence of a lift $\zeta_m$ of $v_n$ 
in $\pXGr\cap J\pXGr$. Since $\pXGr\cap J\pXGr\to M_0$ is a generalized smooth subbundle, it admits a local section $\zeta$ around $m$ in $M_0$ which extends $\zeta_m$.
\end{proof}

\begin{thm}\label{thm:main1}
Let $M$ be a manifold endowed with an action of a Lie group $G$, and $J$ a $G$-invariant generalized complex structure on $M$. Assume that $M_0$ is a $G$-invariant submanifold of $M$ where the $G$-action is free and proper so that the quotient $M_G=M_0/G$ is a manifold.
%If $$\pXGr\cap J\pXGr
%\text{is a generalized smooth subbundle}$$ and
%$$ J\pGXr \cap \pXGr \subset\GXr, $$
Then there is an induced generalized complex structure on $M_G$ if the hypotheses \eqref{5bis} and \eqref{4} of Lemma~\ref{prop:wsmooth} are satisfied.
\end{thm}

Theorem~\ref{thm:main1} has several important consequences.

\begin{cor}
\label{cor:MW}
Under the same hypothesis as in Theorem~\ref{thm:main1}, if moreover
\begin{equation} 
J\pGXr = \GXr ,\label{eqcor1} 
\end{equation} 
then $M_G$ admits an inherited generalized complex structure.
\end{cor}
\begin{proof}
Since $J$ is orthogonal and $\BB=\big(\Bp\big)^{\perp_{\ip{\cdot}{\cdot}}}$,
Eq.~\eqref{eqcor1} implies that $J\pXGr=\XGr$. 
Hence $(\XGr\cap J\pXGr)=\XGr$ is smooth. 
And Eq.~\eqref{4} is satisfied since $\Bp\subset\BB$. The result follows from Theorem~\ref{thm:main1}.
\end{proof}

\begin{rmk}\label{rmk:MW}
This is exactly what happens in the case of Marsden-Weinstein reduction \cite{MW}. In this case,  
$J=\left( \begin{smallmatrix} 0 & -\omega_{\bemol}^{-1} \\ \omega_{\bemol} & 0 \end{smallmatrix}\right)$, 
and $M_0={\mu}^{-1}(0)$, where $\omega\in\Omega^2(M)$ is a symplectic form and $\mu:M\to\mfg^*$
is an equivariant momentum map:
$\hat{A}\ii\omega=d\ip{\mu}{A},\;\forall A\in\mfg$. 
One checks easily that $J(T\orb)=(T M_0)^0$, and therefore
$J(T M_0)^0=T\orb$ since $J^2=I$. Hence the conditions in
Corollary~\ref{cor:MW} are satisfied.
\end{rmk}

More generally one sees that these conditions hold when
$J$ is a symplectic structure twisted by a $B$-field.

\begin{cor}
Let $(M,\omega)$ be a symplectic manifold endowed with an action of a Lie group $G$ preserving the symplectic structure. Assume that there exists a $G$-invariant closed 2-form $B$ on $M$ and an equivariant momentum map $\mu:M\to\mathfrak{g}^*$.
\begin{equation}\label{etoile}\hat{A}\ii\omega=d\ip{\mu}{A}, \quad \forall A\in\mathfrak{g}.
\end{equation}
Assume that $0$ is a regular value of $\mu$ and the $G$-action is free and proper.
Thus $M_0= \mu^{-1}(0)$ is a $G$-invariant submanifold of $M$ and the quotient $M_G=M_0/G$ is a smooth manifold.
 If $i^*B$, where $i: M_0\to M$ is the inclusion,
 is basic with respect to  the $G$-action,
 then the generalized complex structure
 $$J=\begin{pmatrix} 1 & 0 \\ B_{\bemol} & 1 \end{pmatrix}
\begin{pmatrix} 0 & -\omega_{\bemol}^{-1} \\ \omega_{\bemol} & 0 \end{pmatrix}
\begin{pmatrix} 1 & 0 \\ -B_{\bemol} & 1 \end{pmatrix}$$
on $M$ induces a generalized complex structure on $M_G$.
\end{cor}
\begin{proof}
Note that the assumption $i^*B$ being basic implies that
$\BB$ is stable under the bundle map 
$\left( \begin{smallmatrix} 1 & 0 \\ \pm B_{\bemol} & 1 \end{smallmatrix} \right)$ 
of $TM\oplus TM^*$. By the discussion in Remark \ref{rmk:MW}, we know that
$\BB$ is also stable under 
 $\left( \begin{smallmatrix} 0 & -\omega_{\bemol}^{-1} \\ \omega_{\bemol} & 0 \end{smallmatrix} \right)$. The conclusion thus follows 
from Corollary \ref{cor:MW}.
% The bundle endomorphism $J$ is a generalized complex structure \cite[page 45]{Gualtieri}. Since $\omega$ is non degenerate, given $\hat{A}+f\in\GXr$, there exists a unique $z\in TM$ such that $z\ii\omega=\hat{A}\ii B-f$. The 2-form $i^*B$ being basic, \eqref{etoile} implies that $z\in T\orb$ and $\hat{A}\ii\omega + z\ii B\in (TM_0)^0$. Hence $$\begin{pmatrix} 0 & -\omega_{\bemol}^{-1} \\ \omega_{\bemol} & 0 \end{pmatrix}
%\begin{pmatrix} 1 & 0 \\ -B_{\bemol} & 1 \end{pmatrix} \begin{pmatrix} \hat{A} \\ f \end{pmatrix} \in \begin{pmatrix} 1 & 0 \\ -B_{\bemol} & 1 \end{pmatrix} \begin{pmatrix} T\orb \\ (TM_0)^0 \end{pmatrix}$$ and \eqref{4} is satisfied because of the inclusion $\GXr\subset\XGr$.
%Moreover, since $$i^*B \text{ is basic} \quad \Longleftrightarrow \quad TM_0\ii B\subset (T\orb)^0$$
%and $$\hat{A}\ii\omega=d\ip{\mu}{A}, \;\forall A\in\mfg \quad\Longleftrightarrow \quad TM_0\ii\omega\subset (T\orb)^0 ,$$ both $\left( \begin{smallmatrix} 1 & 0 \\ B_{\bemol} & 1 \end{smallmatrix} \right)$ and $\left( \begin{smallmatrix} 0 & -\omega_{\bemol}^{-1} \\ \omega_{\bemol} & 0 \end{smallmatrix} \right)$ preserve $\BB$. Therefore, $J\big(\BB\big)=\BB$. Hence $J\pXGr\cap\pXGr$ is smooth.
\end{proof}

Here is yet another important situation.

\begin{cor}
\label{cor:GS}
Under the same hypothesis as in Theorem \ref{thm:main1}, if moreover
 \begin{equation} 
J\pGXr \cap \pXGr = 0, \label{eqcor2} 
\end{equation}
then $M_G$ admits an inherited generalized complex structure.

In particular, if $J=\left( \begin{smallmatrix} j & 0 \\ 0 & -j^* \end{smallmatrix} \right)$ is a usual complex structure, the above relation is nothing else than the Guillemin-Sternberg  condition \cite{GS}:
 \begin{equation}
\label{eq:GS}
TM=TM_0\oplus j(T\orb).
\end{equation}
In this case, the induced generalized complex structure is still
a complex structure.
\end{cor}
\begin{proof}
Since Eq. \eqref{eqcor2} implies Eq. \eqref{4},
 Condition (\ref{enumv}) of Proposition \ref{prop:condition} holds.
Together  with Eq. \eqref{eqcor2},
this implies that
$$J\pXGr\subset\pXGr \oplus J\pGXr .$$
Since $\Bp\subset\BB$, the projection of $J\pXGr$ onto 
$\XGr$ parallel to $J\pGXr$ is equal to $\big(\BB\big)\cap J\big(\BB\big)$.
Therefore, we have
$$J\pXGr=\big(\BB\big)\cap J\big(\BB\big) \oplus J\pGXr ,$$ 
which implies that $\big(\BB\big)\cap J\big(\BB\big)$ must
be a   smooth subbundle. Thus the conclusion follows 
from Theorem \ref{thm:main1}.

The second part of Corollary \ref{cor:GS} can be easily checked and
is left for the reader.
\end{proof}

\begin{cor}
\label{cor:rie}
Under the same hypothesis as in Theorem \ref{thm:main1}, if moreover
there exists a Riemann metric $g$ such that
${(T\orb\oplus(TM_0)^{\perp})}^\perp\oplus {(T\orb\oplus(TM_0)^{\perp})}^0$
is $J$-stable,  then $M_G$ admits an inherited generalized complex structure.
Here $\perp$ refers to the orthogonal subspace in $TM$ with respect to the metric
$g$.

In particular, if  $(M,g,j)$ is  a $G$-invariant hermitian manifold 
such that  $T\orb\oplus (TM_0)^{\perp}$ is $j$-stable,
 then $M_G$ admits an inherited complex structure.
\end{cor}

\begin{proof}
Since $L:=\big(T\orb\oplus(TM_0)^{\perp} \big)^{\perp}\oplus\big(T\orb\oplus(TM_0)^{\perp} \big)^0$ is $J$-stable, one has $L\subset \pXGr\cap J\pXGr$. Let $Q$ be the projection of 
$(TM\oplus T^*M)_{|M_0}$ onto $(T\orb)^{\perp}\oplus {(TM_0)^{\perp}}^0$ parallel to $\Bp$.
Then, since $\Bp\subset\BB$, one has $$Q\pXGr\subset \big(T\orb\oplus(TM_0)^{\perp} \big)^{\perp}\oplus\big(T\orb\oplus(TM_0)^{\perp} \big)^0=L .$$
Now fix an element $v$ of $TM_G\oplus T^*M_G$, and let $\tilde{v}\in \BB$ be a lift of $v$. Then $Q(\tilde{v})$ is a lift of $v$ in $L$. The bundle $L$ being smooth, there exists a section $s$ of $L$ extending $Q(\tilde{v})$. Hence the hypothesis of Proposition \ref{prop:smoothness} is satisfied and the result follows from Theorem \ref{thm:main}.

In particular, if $J=\left( \begin{smallmatrix} j & 0 \\ 0 & -j^* \end{smallmatrix} \right)$, where $j$ preserves $T\orb\oplus (TM_0)^{\perp}$, and $J$ is orthogonal with respect to $g$, then the smooth subbundle $L$ of $\BB$ is indeed $J$-invariant.
\end{proof}

\begin{rmk}
Note that the induced  complex structure  on the quotient space 
in   K\"ahler reduction can be considered as a special
case of the above  situation.
\end{rmk}

\begin{cor}
Let $(M,j)$ be a complex manifold endowed with an action of a Lie group $G$
 such that $j$ is $G$-invariant. Let $M_0$ be a $G$-invariant submanifold
 of $M$, where the $G$-action is free and proper
 so that the quotient $M_G=M_0/G$ is a manifold. 
Assume that, for any $X_n\in T_n M_G,\;\xi_n\in T^*M_G$,
 there exists a point $m\in M_0$ with $\pi(m)=n$
 and local sections $\tilde{X}\in\secloc{TM_0\cap JTM_0},\;\tilde{\xi}
\in\secloc{(T\orb+j T\orb)^0}$ such that $\pi_*(\tilde{X}|_m)=X_n$
 and $\pi^* \xi_n=\tilde{\xi}|_m$.
Then $M_G$ admits an inherited generalized  complex structure, which
is still a  complex structure.

In particular, the above conditions are satisfied if
$TM_0\cap jTM_0$ and $T\orb+ j T\orb$
are generalized smooth subbundles  of  $TM_0$.
\end{cor}
\begin{proof}
The above hypothesis is the projection of the hypothesis of Proposition
 \ref{prop:smoothness} onto the tangent and cotangent bundles. 
\end{proof}

\section{Description of $J_G$}
\label{sec:jg}
This section is devoted to an explicit description of the
induced generalized complex structure
$J_G$ under the assumption that the hypothesis of Theorem \ref{thm:main}
holds.

Let $B$ denote the vector bundle $\BB$ and $B^{\perp}$
 its $\ip{\cdot}{\cdot}$-orthogonal subbundle: $\Bp$.
As in Corollary \ref{cor:directsum}, we know that every vector in
 $TM_G\oplus T^*M_G$ has a lift in $B\cap JB$.
In other words, the restriction of the bundle map $\Pi$ (as in Eq. \eqref{eq:liftbundlemap}) to $B\cap JB$ is surjective. 
Then $\Pi_{|B\cap JB}$ induces a $G$-invariant vector bundle map 
$$\xymatrix{ \tfrac{B\cap JB}{\ker(\Pi_{|B\cap JB})} \ar[d] \ar[r]^{\Xi} & TM_G\oplus T^*M_G \ar[d] \\ 
M_0 \ar[r]^{\pi} & M_G }$$ 
which is indeed a pullback bundle. Note that $\ker_{|B\cap JB}=B^{\perp}\cap JB$. 
Now, by Corollary \ref{cor:directsum} and Proposition \ref{prop:condition},
 we have $JB^{\perp}\cap B\subset B^{\perp}$. Moreover
 since $B^{\perp}\subset B$, we thus have $JB^{\perp}\cap B\subset B^{\perp}\cap JB$. Since
 $J^2=-I$, it follows that  $JB^{\perp}\cap B=B^{\perp}\cap JB$.
It is simple to see that  $B^{\perp}\cap JB$ is $J$-stable.
 Therefore, $J_{|B\cap JB}$ induces a bundle isomorphism 
$$\xymatrix{ \tfrac{B\cap JB}{B^{\perp}\cap JB} \ar[d] \ar[r]^{J'} & \tfrac{B\cap JB}{B^{\perp}\cap JB} \ar[d] \\ 
M_0 \ar[r]^{id} & M_0 }$$ 
whose $G$-invariance implies that it factorizes through $\Xi$ so
that it induces a  bundle map $J_G$ as the inherited generalized
 complex structure 
$$\xymatrix{ TM_G\oplus T^*M_G \ar[d] \ar[r]^{J_G} & TM_G\oplus T^*M_G \ar[d] \\ M_G \ar[r]^{id} & M_G .}$$

Thus we obtain the following

\begin{thm}
Under the same hypothesis of Theorem \ref{thm:main}, the
induced generalized complex structure $J_G$ on $M_G$ can
be described by the  following commutative diagram. 
$$\xymatrix{ TM\oplus T^*M \ar[r]^J & TM\oplus T^*M \\ 
B\cap JB \ar[u] \ar[d] \ar[r]^J & B\cap JB \ar[u] \ar[d] \\ 
\tfrac{B\cap JB}{B^{\perp}\cap JB}\ar[d]_{\Xi} \ar[r]^{J'} & \tfrac{B\cap JB}{B^{\perp}\cap JB}\ar[d]^{\Xi} \\ 
TM_G\oplus T^*M_G \ar[r]^{J_G} & TM_G\oplus T^*M_G }$$
Hence, if $\tilde{v}\in B\cap JB$ is a lift of $v\in TM_G\oplus T^*M_G$, then 
$J\tilde{v}\in B\cap JB$ is a lift of $J_G v\in TM_G\oplus T^*M_G$.
\end{thm}

\begin{rmk}
In \cite{math.DG/0412097}, Crainic showed that a generalized complex structure
 $J$ is an endomorphism of the generalized tangent
 bundle $TM\oplus T^*M$ built out of a Poisson bivector $\pi$, an endomorphism $N$ of the tangent bundle and a 2-form $\sigma$:
 $$J=\begin{pmatrix} N & \pi^{\sharp} \\ \sigma_{\flat} & -N^* \end{pmatrix}$$ with $N^2+\pi^{\sharp} \sigma_{\flat }=-\idn$, 
which satisfy some compatibility conditions 
 resembling those given by Kosmann-Schwarzbach, Magri and Morosi
\cite{MR1077465, MR773513}  in their definition of Poisson Nijenhuis structures
(see also \cite{MR2163575, SX}). It would be interesting to  investigate
the meaning of each component for the reduction $J_G$.
\end{rmk}

\section{Reduction of Generalized K\"ahler structures}

\emph{Generalized K\"ahler structures} as introduced by
Hitchin \cite{Hitchin} and Gualtieri \cite{Gualtieri}
consists of  a   pair of commuting generalized complex structures
 $(J_1,J_2)$ such that $< J_1J_2\; \cdot , \cdot>$ is  positive definite.

\begin{example}
If $(M,\omega,j)$ is a usual K\"ahler manifold,
 then $$J_1=\begin{pmatrix} 0 & -\omega_{\bemol}^{-1} \\ \omega_{\bemol} & 0 
\end{pmatrix} \qquad J_2=\begin{pmatrix} j & 0 \\ 0 & -j^* \end{pmatrix}$$ defines a generalized K\"ahler structure.
\end{example}

As in Section \ref{sec:rds}, we assume that 
  $M_0$ is a $G$-invariant submanifold of $M$ where the $G$-action is
free and proper so that the quotient $M_G=M_0/G$ is a manifold.
Since $TM\oplus T^*M$ decomposes as  the direct sum of four
different subbundles with  each being
 constituted of the common eigenvectors of $J_1$ and $J_2$ determined by a given pair of eigenvalues, it is natural to ask if the reduction procedure outlined in Section \ref{sec:rds} can be applied to $J_1$ and $J_2$ simultaneously to get a generalized K\"ahler structure on the quotient $M_G$.

Let us introduce some notations.
By $E^{\pm}$ and $E_{\pm}$ we denote the $\pm i$-eigenbundles of $J_1$
and $J_2$ respectively. Let $\BC=B\otimes\mathbb{C}$. For any subbundle $V$ of $\TC M\oplus \TC^*M$, $\Phi(V)$  denotes
 the subbundle $\Pi(V\cap \BC)$ of $\TC M_G\oplus \TC^*M_G$,
where $\Pi$ is the bundle map as in Eq. \eqref{eq:liftbundlemap}.

\begin{prop}
\label{pro:J12}
Assume that $(J_1, J_2)$ are a $G$-invariant 
 K\"ahler structure. 
\begin{enumerate}
\item The following two statements are equivalent: 
\begin{equation}
\label{eq:pm}
\TC M_G\oplus \TC^*M_G=\Phi(E^+\cap E_+)+\Phi(E^+\cap E_-)+\Phi(E^-\cap E_+)+\Phi(E^-\cap E_-)
\end{equation}
and 
$$ \Pi(\zambon)=\ctXG, $$
where, as in Section \ref{sec:jg},  $B=\BB$.
\item In case the above condition is satisfied, then the above sum is direct and $$\Phi(E^{\pm}\cap E_{\pm})=\Phi(E^{\pm})\cap\Phi(E_{\pm}) .$$
\end{enumerate}
\end{prop}

\begin{proof}
\begin{enumerate}
\item 
\begin{itemize}
\item[$\boxed{\Rightarrow}$] 
Take $v\in\ctXG$. By assumption, it can be decomposed 
as $$v=\Pi(\tilde{v}^+_+)+\Pi(\tilde{v}^+_-)+\Pi(\tilde{v}^-_+)+\Pi(\tilde{v}^-_-)$$
where $\tilde{v}^{\pm}_{\pm}\in E^{\pm}\cap E_{\pm}\cap \BC$. 
Let $$\tilde{v}=
\tilde{v}^+_++\tilde{v}^+_-+\tilde{v}^-_++\tilde{v}^-_-$$
Then $\tilde{v}\in \BC$ and
  $\Pi(\tilde{v})=v$. Hence $\tilde{v}$ is a lift of $v$.
Moreover, 
$$J_1\tilde{v}=i(\tilde{v}^+_++\tilde{v}^+_--\tilde{v}^-_+-\tilde{v}^-_-)\in \BC \quad \text{and} \quad 
J_2\tilde{v}=i(\tilde{v}^+_+-\tilde{v}^+_-+\tilde{v}^-_+-\tilde{v}^-_-)\in \BC .$$
Hence, $\tilde{v}\in J_1\BC$ and $\tilde{v}\in J_2 \BC$
 since $J_1^2=J_2^2=-I$.
Therefore $\REAL(\tilde{v})\in\zambon$  is a lift of $v\in\ctXG$.
\item[$\boxed{\Leftarrow}$] 
Conversely, if any $v\in \ctXG$ has a lift $\tilde{v}$ in $\zambon$, one has 
\begin{multline*} v=\Pi\big(\tfrac{1}{2}(I-iJ_1)\tfrac{1}{2}(I-iJ_2)\tilde{v}\big)
+\big(\tfrac{1}{2}(I-iJ_1)\tfrac{1}{2}(I+iJ_2)\tilde{v}\big) \\
+\big(\tfrac{1}{2}(I+iJ_1)\tfrac{1}{2}(I-iJ_2)\tilde{v}\big)
+\big(\tfrac{1}{2}(I+iJ_1)\tfrac{1}{2}(I+iJ_2)\tilde{v}\big) .
\end{multline*}
\end{itemize}
Eq. \eqref{eq:pm} thus follows.

\item 
Finally, since 
\begin{equation} \label{a} \Phi(E^{\pm}\cap E_{\pm})\subseteq\Phi(E^{\pm})\cap\Phi(E_{\pm}) ,\end{equation}
the condition 
\begin{equation} \label{b} \cctXG=\Phi(E^+\cap E_+)+\Phi(E^+\cap E_-)+\Phi(E^-\cap E_+)+\Phi(E^-\cap E_-) \end{equation} 
implies that 
\begin{equation} \label{c} \cctXG=\Phi(E^+)\cap \Phi(E_+)+\Phi(E^+)\cap \Phi(E_-)+\Phi(E^-)\cap \Phi(E_+)+\Phi(E^-)\cap \Phi(E_-) \end{equation}
and 
\begin{align*} \cctXG&=\Phi(E^+)+\Phi(E^-), & \cctXG&=\Phi(E_+)+\Phi(E_-) .\end{align*}
Applying Proposition~\ref{prop:condition} to $J_1$ and $J_2$, we obtain 
$$ \Phi(E^+)\cap\Phi(E^-)=0 \qquad \text{and } \qquad \Phi(E_+)\cap\Phi(E_-)=0 .$$
This implies that the sums in \eqref{b} and \eqref{c} must be direct. 
Hence Eq.~\eqref{a} implies that $\Phi(E^{\pm}\cap E_{\pm})
=\Phi(E^{\pm})\cap\Phi(E_{\pm})$.
\end{enumerate}
\end{proof}

\begin{thm}\label{thm:6.2}
Let $(M,J_1,J_2)$ be  a $G$-invariant generalized  K\"ahler manifold. Assume that $M_0$ is a $G$-invariant submanifold of $M$ where the $G$-action is free and proper so that the quotient $M_G=M_0/G$ is a manifold.
If, for any  $v_n\in\ctXG$, there exists a point $m\in M_0$ and a local smooth section $\zeta$ of $\zambon\to M_0$ around $m$ such that $\Pi(\zeta|_m)=v_n$, then there exists an inherited generalized K\"ahler structure $({J_1}_G,{J_2}_G)$ on $M_G$ such that 
\begin{equation}\label{eq:JG12} 
\begin{gathered} 
\Pi(J_1 w)= {J_1}_G \Pi(w) \\ \Pi(J_2 w)= {J_2}_G \Pi(w) \end{gathered} \qquad \forall w\in\zambon .
\end{equation} 

In particular, this condition holds if
\begin{itemize}
\item $B=\zambon+B^{\perp}$ and
\item $\zambon$ is a generalized smooth subbundle.
\end{itemize}
\end{thm}

\begin{proof}
According to  Theorem \ref{thm:main}, the generalized 
complex structures $J_1$ and $J_2$ induce a pair of generalized complex
 structures ${J_1}_G$ and ${J_2}_G$ on $M_G$ satisfying Eq. \eqref{eq:JG12}.
  Let $E'^{\pm}$ and $E'_\pm$ denote their 
$\pm \imaginary$-eigenbundles respectively. By Proposition \ref{pro:J12}, we have 
$$\TC M_G\oplus \TC^* M_G =\big({E'}^{+}\cap E'_+\big)\oplus \big({E'}^+\cap E'_-\big)\oplus \big({E'}^-\cap E'_+\big)\oplus \big({E'}^-\cap E'_-\big),$$
which implies that ${J_1}_G$ and ${J_2}_G$ commute.

Finally, for $v\in \ctXG$, choose $w\in\zambon$. By Eq.~\eqref{eq:JG12}, we have
$$<{J_1}_G{J_2}_G v,v>=<{J_1}_G{J_2}_G\Pi(w),\Pi(w)>
=<\Pi(J_1 J_2 w),\Pi(w)>=<J_1 J_2 w,w>,$$
where the last equality follows from the fact that
$\Pi|_B$ preserves the scalar products.
It thus follows that $<{J_1}_G {J_2}_G \;\cdot,\cdot>$ is
positive definite.
This concludes the proof.
\end{proof}

The following corollary is to Theorem~\ref{thm:6.2} what Corollary~\ref{cor:rie} is to theorem~\ref{thm:main}.

\begin{cor}
\label{cor:rie1}
Let $(M,J_1,J_2)$ be a $G$-invariant generalized K\"ahler manifold.
Assume that $M_0$ is a $G$-invariant submanifold of $M$,
 where the $G$-action is free and proper so that the quotient $M_G=M_0/G$
 is a manifold. Assume that  there exists a Riemann metric $g$ such that
${(T\orb\oplus(TM_0)^{\perp})}^\perp\oplus {(T\orb\oplus(TM_0)^{\perp})}^0$
is stable under both $J_1$ and $J_2$,  then $M_G$ admits an inherited generalized complex structure.
\end{cor}

\begin{proof}
The proof is exactly the same as for Corollary~\ref{cor:rie} except that, here, since $L$ is both $J_1$- and $J_2$-invariant, one has $L\subset\zambon$. Then the conclusion follows from Theorem~\ref{thm:6.2}.
\end{proof}

As proved by Gualtieri \cite{Gualtieri}, a generalized K\"ahler structure $(J_1,J_2)$ on a manifold $M$
is fully characterized by a quadruple $(g, b, J_+, J_-)$, where $g$ is a Riemannian
metric, $b$ is a two form and $J_{\pm}$ are two (integrable) complex structures on $M$, compatible with $g$, such that 
$$db(X,Y,Z)=d\omega_{\pm}(J_{\pm} X,J_{\pm} Y,J_{\pm} Z), \quad \forall X,Y,Z\in\mathfrak{X}(M)$$ with $\omega_{\pm}\in \Omega^2(M)$ defined by $\ip{\omega_{\pm}X}{Y}=\pm g(J_{\pm}X,Y)$. These data explicitly determine the pair of generalized complex structures: 
\begin{gather*}
J_1=\tfrac{1}{2}\begin{pmatrix}1 & 0 \\ b & 1 \end{pmatrix}
\begin{pmatrix} J_++J_- & -(\omega_+^{-1}+\omega_-^{-1}) \\ \omega_++\omega_- & -(J_+^*+J_-^*) \end{pmatrix}
\begin{pmatrix}1 & 0 \\ -b & 1 \end{pmatrix} \\
J_2=\tfrac{1}{2}\begin{pmatrix}1 & 0 \\ b & 1 \end{pmatrix}
\begin{pmatrix} J_+-J_- & -(\omega_+^{-1}-\omega_-^{-1}) \\ \omega_+-\omega_- & -(J_+^*-J_-^*) \end{pmatrix}
\begin{pmatrix}1 & 0 \\ -b & 1 \end{pmatrix} .
\end{gather*}
Here, for simplicity, we identify a 2-form with its associated bundle map.

Using exactly this same $g$ and applying Corollary~\ref{cor:rie1}, we are led to the main result of this section,
which generalizes the usual K\"ahler reduction \cite{GS, Hitchinet.al., K}.

\begin{thm} 
Let $(M, J_1, J_2)$ be  a $G$-invariant generalized  K\"ahler manifold.
Assume that $M_0$ is a $G$-invariant submanifold of $M$, where the $G$-action is free and proper so that the quotient $M_G=M_0/G$ is a manifold.
Let $(g,b,\omega_+,\omega_-)$ be the quadruple associated to a generalized K\"ahler structure $(J_1,J_2)$ as in
\cite{Gualtieri}. If $\omega_{\pm}(T\orb)=(TM_0)^0$
and $(T\orb)^{\perp}\cap TM_0$ is stable under $\omega_{\pm}^{-1}b$,
then $M_G$ inherits a generalized K\"ahler structure.
\end{thm}

\begin{proof}
It is sufficient to check that the conditions in Corollary~\ref{cor:rie1} are satisfied.

First, since $\omega_{\pm}$ are antisymmetric and non-degenerate, $\omega_{\pm}(T\orb)=(TM_0)^0$ implies  $\omega_{\pm}(TM_0)=(T\orb)^0$.

Second, $\omega_{\pm}(T\orb)=(TM_0)^0$ also implies $\omega_{\pm}(T\orb^{\perp})=\big(TM_0^{\perp}\big)^0$ because $J_{\pm}$ are isometric with respect to $g$. Therefore, 
\begin{equation} \omega_{\pm}\big(T\orb\oplus (TM_0)^{\perp}\big)^{\perp}=\omega_{\pm}(T\orb^{\perp}\cap TM_0)=\big(TM_0^{\perp}\big)^0\cap \big(T\orb\big)^0=\big(T\orb\oplus (TM_0)^{\perp}\big)^0 \label{intersection} .\end{equation}

Finally, one easily checks that the image of any $$v+f\in \big(T\orb\oplus (TM_0)^{\perp}\big)^{\perp}\oplus \big(T\orb\oplus (TM_0)^{\perp}\big)^0$$ under $$J_1\pm J_2=\begin{pmatrix}1 & 0 \\ b & 1 \end{pmatrix}
\begin{pmatrix} J_{\pm} & -\omega_{\pm}^{-1} \\ \omega_{\pm} & -J_{\pm}^* \end{pmatrix}
\begin{pmatrix}1 & 0 \\ -b & 1 \end{pmatrix}$$ lies in $\big(T\orb\oplus (TM_0)^{\perp}\big)^{\perp}\oplus \big(T\orb\oplus (TM_0)^{\perp}\big)^0$ if and 
only if both Eq.~\eqref{intersection} and the condition $\omega_{\pm}^{-1}b((T\orb)^{\perp}\cap TM_0)\subset (T\orb)^{\perp}\cap TM_0$ hold.
\end{proof}

Note that $\omega_{\pm}(T\orb) =(TM_0)^0$ is satisfied when
$M_0$ is the zero level set of an equivariant momentum map.
Thus the theorem above reduces to the usual K\"ahler reduction when $M$ is a K\"ahler manifold \cite{GS, Hitchinet.al., K}.

%\bibliography{redbib2}
%\bibliographystyle{plain}

\end{document}